\documentclass[12pt]{amsproc}
\usepackage{amssymb}
\usepackage{amsmath}
\usepackage[english,russian]{babel}

\begin{document}



{\Large \begin{center} Odinary differential operators of odd order with distribution coefficients
 \end{center}}

\bigskip
\centerline{K.~A.Mirzoev and A.~A.~Shkalikov}\footnote{This work is supported by Russian Science Foundation (RNF) under grant No 17-11-01215.}

\bigskip
\medskip

\centerline{{\bf Abstract}}
\bigskip
We work with differential expressions of the form
\begin{multline*}
\tau_{2n+1} y=(-1)^ni
\{(q_{0}y^{(n+1)})^{(n)}+(q_{0}y^{(n)})^{(n+1)}\}+
\sum\limits_{k=0}^{n}(-1)^{n+k}(p^{(k)}_ky^{(n-k)})^{(n-k)}+\\
+i\sum\limits_{k=1}^{n}(-1)^{n+k+1}\{(q^{(k)}_{k}y^{(n+1-k)})^{(n-k)}+
(q^{(k)}_{k}y^{(n-k)})^{(n+1-k)}\},
\end{multline*}
where the complex valued coefficients $p_j$ and $q_j$ are subject  the following conditions:  $ q_0(x) \in AC_{loc}(a,b)$, $Re \,q_0>0$, while all the other functions
 $$
 q_1(x),q_2(x),\ldots,q_{n}(x), p_0(x),p_1(x),\ldots,p_n(x)
 $$
 belong to the space $L^1_{loc}(a,b)$. This implies that the coefficients $p^{(k)}_{k}$ and $q^{(k)}_{k}$ in the expression $\tau_{2n+1}$ are  distributions of singularity order $k$.

  The main objective of the paper  is to represent the differential expression $\tau_{2n+1}$ in the other (regularized) form which allows to define the minimal and maximal operators associated with this differential expression.


\bigskip
\medskip

Цель этой заметки  --- научиться определять дифференциальные операторы  нечетного порядка с коэффициентами-распределениями при некоторых
естественных условиях на эти коэффициенты. Мы будем действовать по той же схеме, что в нашей прежней работе \cite{MSh},  в которой аналогичная задача
решалась для операторов четного порядка. Одновременно отметим, что в нечетном случае имеется своя специфика и новые трудности.

Рассмотрим  дифференциальное выражение нечетного порядка
\begin{multline}\label{0}
\tau_{2n+1} y=(-1)^ni
\{(q_{0}y^{(n+1)})^{(n)}+(q_{0}y^{(n)})^{(n+1)}\}+
\sum\limits_{k=0}^{n}(-1)^{n+k}(p^{(k)}_ky^{(n-k)})^{(n-k)}+\\
+i\sum\limits_{k=1}^{n}(-1)^{n+k+1}\{(q^{(k)}_{k}y^{(n+1-k)})^{(n-k)}+
(q^{(k)}_{k}y^{(n-k)})^{(n+1-k)}\},
\end{multline}
предполагая, что комплекснозначные  коэффициенты $p_j$ и $q_j$  таковы, что  функция
 $ q_0(x) \in AC_{loc}(a,b)$, $Re \,q_0>0$, а все функции
 $$ q_1(x),q_2(x),\ldots,q_{n}(x), p_0(x),p_1(x),\ldots,p_n(x)
 $$
 принадлежат пространству $L^1_{loc}(a,b)$. Тем самым, функции $p^{(k)}_{k}$ и $q^{(k)}_{k}$, участвующие в качестве коэффициентов в \eqref{0}, являются распределениями порядка сингулярности $k$.

Введем обозначения
 $$
  \varphi_k:=p_k+iq_{k},\qquad  \tilde{\varphi_k}:=p_k-iq_{k},\qquad k=1,2,\ldots,n,
  $$
 и определим квадратную
матрицу $F_{2n+1}$, полагая
$$ \small{\setcounter{MaxMatrixCols}{15} {F}_{2n+1}=\begin{pmatrix} 0 & 1&  .& .& .& 0&0&0&0&.&.&
.&0\\ 0& 0 &   .& .& .&0&0&0 &0&.& .&.&0  \\ .&. &  . &. & .&.&.&. &.&.&.&.&.\\ 0&0 &
.& .& .& 1&0&0 &0&.& .&.&0\\  0&0 &
.& .& .& 0&f_{n,n+1}&0&0 &.& .&.&0\\  f_{n+1,1} & f_{n+1,2}  &.&.&.&f_{n+1,n}  & f_{n+1,n+1}
&  f_{n+1,n+2} & 0 &.&.&.&0\\ 0 & f_{n+2,2}   &.&.& .&f_{n+2,n}  & f_{n+2,n+1}
& 0 &1&.&.&.&0\\ .&. &. & . & . &.&.&. &.&.&. &.\\ 0&0
&.& .&.&f_{2n,n} & f_{2n,n+1} & 0 &0&.&.&.&1\\ 0& 0
&.&.& .&0 & f_{2n+1,n+1} & 0 &0&.&.&.&0\\ \end{pmatrix} }, $$
где ненулевые элементы $f_{ij}$ этой матрицы  определяются равенствами
$$
f_{i,i+1}=1, \qquad  i=1,2, \ldots,n-1,n+2, \ldots,2n,
$$
$$
f_{n,n+1}=f_{n+1,n+2}=\frac{1}{\sqrt{2q_0}},\qquad f_{n+1,n+1-j}=(-1)^j\frac{i\tilde{\varphi}_j}{\sqrt{2q_0}}, \, j=1,2, \ldots n, $$ $$ f_{n+1,n+1}=\frac{ip_0}{2q_0},$$
$$ f_{n+1+j,n+1}= \frac{i \varphi_j}{\sqrt{2q_0}},  j=1,2, \ldots n,$$
$$f_{n+1+k,n+k-j}=(-1)^{j+k+1}C_{j+1}^{k}\left[ip_{j+1}+(1-\frac{2k}{j+1})q_{j+1}\right], $$
 $$k=1,2,\ldots,n-1,\, j=k,k+1,\ldots,n-1.
$$

 Для наглядности приведем явный вид матриц $F_3, F_5$ и $F_7$:

$$ \large F_3=\begin{pmatrix}
0&\frac{1}{\sqrt{2q_0}}&0\\
-\frac{i\tilde{\varphi_1}}{\sqrt{2q_0}}&\frac{ip_0}{2q_0}&\frac{1}{\sqrt{2q_0}}\\
0&\frac{i\varphi_1}{\sqrt{2q_0}}&0\\\end{pmatrix}, \quad
F_5=\begin{pmatrix}
0&1&0  \quad&0\quad&0\\
0&0&\frac{1}{\sqrt{2q_0}}\quad&0\quad&0\\
\frac{i\tilde{\varphi}_2}{\sqrt{2q_0}}&-\frac{i\tilde{\varphi}_1}{\sqrt{2q_0}}&\frac{ip_0}{2q_0} \quad&\frac{1}{\sqrt{2q_0}}\quad &0\\
0&-2ip_2&\frac{i{\varphi}_1}{\sqrt{2q_0}} \quad &0\quad&1\\
0&0&\frac{i{\varphi}_2}{\sqrt{2q_0}}\quad &0\quad&0\\
\end{pmatrix}
$$

и

$$
F_7=\begin{pmatrix}
0\qquad&1\qquad&0 \qquad &0\qquad&0\qquad&0\qquad&0\\
0\qquad&0\qquad&1\qquad  &0\qquad&0\qquad&0\qquad&0\\
0\qquad&0\qquad&0\qquad&\frac{1}{\sqrt{2q_0}}\qquad&0\qquad&0\qquad&0\\
-\frac{i\tilde{\varphi}_3}{\sqrt{2q_0}}\qquad&\frac{i\tilde{\varphi}_2}{\sqrt{2q_0}}\qquad&-\frac{i\tilde{\varphi}_1}{\sqrt{2q_0}}\qquad&\frac{ip_0}{2q_0} \qquad &\frac{1}{\sqrt{2q_0}}\qquad &0\qquad&0\\
0\qquad&q_3+3ip_3\qquad&-2ip_2\qquad&\frac{i{\varphi}_1}{\sqrt{2q_0}}\qquad  &0\qquad&1\qquad&0\\
0\qquad&0\qquad&q_3-3ip_3\qquad&\frac{i{\varphi}_2}{\sqrt{2q_0}} \qquad&0\qquad&0\qquad&1\\
0\qquad&0\qquad&0\qquad&\frac{i{\varphi}_3}{\sqrt{2q_0}}\qquad &0\qquad&0\qquad&0\\
\end{pmatrix}.
$$

Матрица $F_{2n+1}$ имеет специальный вид. Такого типа матрицы рассматривались и ранее, и в литературе их называют матрицами Шина--Зетла (см., например, \cite{EM}). В работе  \cite{MSh}  класс таких матриц  обозначался через $ \mathcal{S}_{2n+1} (a,b)$.

 Определим теперь квазипроизводные $y^{[j]}$, $j=0,1,...2n ,$
 и квазидифференциальное выраженние $\tau y$ для заданной функции $y$, согласованные с матрицей $F_{2n+1}$.  Элементарные вычисления показывают, что
\begin{equation} \label{f00} y^{[k]}=y^{(k)},\,\,\,\,\,\,\,\,\,k=0,1,\ldots,n-1, \end{equation}
\begin{equation} \label{f10} y^{[n]}=\sqrt{2q_0}y^{(n)},
\end{equation}
\begin{equation} \label{f201}
y^{[n+1]}=
\sqrt{2q_0}(\sqrt{2q_0}y^{(n)})^{\prime}-ip_0y^{(n)}-i\sum\limits_{j=1}^{n}(-1)^j\tilde{\varphi_j}y^{(n-j)},
\end{equation}
\begin{equation} \label{f20}
y^{[n+k+1]}=
(y^{[n+k]})^{\prime}-\sum\limits_{j=k-1}^{n-1}(-1)^{j+k+1}C_{j+1}^{k}\left[ip_{j+1}+(1-\frac{2k}{j+1})q_{j+1}\right]y^{(n+k-j-1)},
\end{equation}
$$
k=1,2,\ldots,n-1,
$$
и
\begin{equation} \label{f30}
 \tau_{2n+1}
y=i^{2n+1}[(y^{[2n]})^{\prime}-i{\varphi_n}y^{(n)}].
\end{equation}

Сначала приведем к дивергентной форме вида \eqref{0} выражение $\tau_m
y$ при значениях $m=3,5$ и $7$.

Пусть $m=3$, тогда, применяя формулы \eqref{f00}--\eqref{f30}, получим
$$
y^{[1]}=\sqrt{2q_0}y^{\prime},
$$
$$
y^{[2]}=\sqrt{2q_0}(\sqrt{2q_0}y^{\prime})^{\prime} -ip_0y^{\prime}+i\tilde{\varphi_1}y =\{(q_0y^{\prime})^{\prime}+q_0y^{\prime\prime}\}-ip_0y^{\prime}+i\tilde{\varphi_1}y=\Phi_0+i\tilde{\varphi_1}y,
$$
где $$\Phi_0=\{(q_0y^{\prime})^{\prime}+q_0y^{\prime\prime}\}-ip_0y^{\prime},$$
и $$
\tau_3 y=i^3\Bigl((y^{[2]})^{\prime}-i\varphi_1y^{\prime}\Bigr)
=i^3\Bigl((\Phi_0^{\prime}+i(\tilde{\varphi_1}y)^{\prime}-i\varphi_1y^{\prime}\Bigr)=i^3\Bigl((\Phi_0^{\prime}+\{(q_{1}y)^{\prime}+q_{1}y^{\prime}\}+i
p^{\prime}_1y\Bigr).
$$Таким образом,
$$
\tau_3 y=i^3\Bigl(\{(q_{0}y^{\prime})^{\prime\prime}+
(q_{0}y^{\prime\prime})^{\prime}\}-i(p_0y^{\prime})^{\prime}+\{(q_{1}y)^{\prime}+q_{1}y^{\prime}\}+i
p^{\prime}_1y\Bigr)=
$$
$$
=-i\{(q_{0}y^{\prime})^{\prime\prime}+
(q_{0}y^{\prime\prime})^{\prime}\}-(p_0y^{\prime})^{\prime}-i\{(q_{1}y)^{\prime}+q_{1}y^{\prime}\}+
p^{\prime}_1y.
$$

Пусть теперь $m=5$. Согласно формулам \eqref{f00}--\eqref{f30}, имеем
$$
y^{[1]}=y^{\prime},
\qquad
y^{[2]}=\sqrt{2q_0}y^{\prime\prime},
$$
$$
y^{[3]}=\sqrt{2q_0}(\sqrt{2q_0}y^{\prime\prime})^{\prime}-ip_0y^{\prime\prime}+i\tilde{\varphi_1}y^{\prime}-i\tilde{\varphi_2}y
=$$
$$=\{(q_0y^{\prime\prime})^{\prime}+q_0y^{\prime\prime\prime}\}-ip_0y^{\prime\prime}+i\tilde{\varphi_1}y^{\prime}-i\tilde{\varphi_2}y=
\Phi_0+i\tilde{\varphi_1}y^{\prime}-i\tilde{\varphi_2}y,
$$
где $$\Phi_0=\{(q_0y^{\prime\prime})^{\prime}+q_0y^{\prime\prime\prime}\}-ip_0y^{\prime\prime},$$

$$
y^{[4]}=\Phi_0^{\prime}+\{(q_{1}y^{\prime})^{\prime}+q_{1}y^{\prime\prime}\}+i
p^{\prime}_1y^{\prime}-i(\tilde{\varphi_2}y)^{\prime}+2ip_2y^{\prime}=\Phi_1-i(\tilde{\varphi_2}y)^{\prime}+2ip_2y^{\prime},
$$
где $$\Phi_1=\Phi_0^{\prime}+\{(q_{1}y^{\prime})^{\prime}+q_{1}y^{\prime\prime}\}+i
p^{\prime}_1y^{\prime},$$

и поэтому
$$
\tau_5 y=i^5\Bigl((y^{[4]})^{\prime}-i\varphi_2y^{\prime\prime}\Bigr)=
i^5\Bigl( \Phi_1^{\prime}-i(\tilde{\varphi_2}y)^{\prime\prime}+2i(p_2y^{\prime})^{\prime}-i(\tilde{\varphi_2}y)^{\prime\prime}\Bigr)=
$$
$$
=i^5\Bigl( \Phi_1^{\prime}-\{(q_{2}^{\prime}y^{\prime})+
(q_{2}^{\prime}y)^{\prime}\}-ip_2y^{\prime\prime}\Bigr)
$$
Таким образом,
$$
\tau_5 y=i\{(q_{0}y^{\prime\prime})^{\prime\prime\prime}+
(q_{0}y^{\prime\prime\prime})^{\prime\prime}\}+(p_0y^{\prime\prime})^{\prime\prime}+i\{(q_{1}y^{\prime})^{\prime\prime}+(q_{1}y^{\prime\prime})^{\prime}\}-
(p^{\prime}_1y^{\prime})^{\prime}-i\{(q_{2}^{\prime}y^{\prime})+
(q_{2}^{\prime}y)^{\prime}\}+p_2^{\prime\prime}y
$$

Рассмотрим далее случай, когда $m=7$,
$$
y^{[1]}=y^{\prime},
\qquad
y^{[2]}=y^{\prime\prime},
\qquad
y^{[3]}=\sqrt{2q_0}y^{\prime\prime\prime},
$$

$$
y^{[4]}=\sqrt{2q_0}(\sqrt{2q_0}y^{\prime\prime\prime})^{\prime}-ip_0y^{\prime\prime\prime}+
i\tilde{\varphi_3}y-i\tilde{\varphi_2}y^{\prime}+i\tilde{\varphi_1}y^{\prime\prime}=
$$
$$=\{(q_0y^{\prime\prime\prime})^{\prime}+q_0y^{(4)}\}-ip_0y^{\prime\prime\prime}+
i\tilde{\varphi_3}y-i\tilde{\varphi_2}y^{\prime}+i\tilde{\varphi_1}y^{\prime\prime}=
\Phi_0+
i\tilde{\varphi_3}y-i\tilde{\varphi_2}y^{\prime}+i\tilde{\varphi_1}y^{\prime\prime},
$$
где $$\Phi_0=\{(q_0y^{\prime\prime\prime})^{\prime}+q_0y^{(4)}\}-ip_0y^{\prime\prime\prime},$$

$$
y^{[5]}=\Phi_1-i(\tilde{\varphi_2}y^{\prime})^{\prime}+i(\tilde{\varphi_3}y)^{\prime}-(q_3+3ip_3)y^{\prime}+2ip_2y^{\prime\prime},
$$
где $$\Phi_1=\Phi_0^{\prime}+\{(q_{1}y^{\prime\prime})^{\prime}+q_{1}y^{\prime\prime\prime}\}+i
p^{\prime}_1y^{\prime\prime},$$
а
$$
y^{[6]}=\Phi_2+i(\tilde{\varphi_3}y)^{\prime\prime}-((q_3+3ip_3)y^{\prime})^{\prime}-(q_3-3ip_3)y^{\prime\prime},
$$
где $$\Phi_2=\Phi_1^{\prime}-\{(q^{\prime}_{2}y^{\prime})^{\prime}+q^{\prime}_{2}y^{\prime\prime}\}-i
p^{\prime\prime}_2y^{\prime},$$
и
$$
\tau_7 y=i^7\Bigl((y^{[6]})^{\prime}-i\varphi_3y^{\prime\prime\prime}\Bigr)=
i^7\Bigl( \Phi_2^{\prime}+\{(q^{\prime\prime}_{3}y )^{\prime}+q^{\prime\prime}_{3}y^{\prime}\}+ip^{\prime\prime\prime}_3y\Bigr).
$$

Таким образом,
$$
\tau_7 y=i^7\Bigl(\{(q_{0}y^{\prime\prime\prime})^{(4)}+
(q_{0}y^{(4)})^{\prime\prime\prime}\}-i(p_0y^{\prime\prime\prime})^{\prime\prime\prime}+$$
$$
+\{(q_{1}y^{(\prime\prime)})^{\prime\prime\prime}+(q_{1}y^{\prime\prime\prime})^{\prime\prime}\}
+i(p^{\prime}_1y^{\prime\prime})^{\prime\prime}-
\{(q^{\prime}_{2}y^{\prime})^{\prime\prime}+(q^{\prime}_{2}y^{\prime\prime})^{\prime  }\}-i(p^{\prime\prime}_2y^{\prime})^{\prime}+
\{(q^{\prime\prime}_{3}y )^{\prime}+(q^{\prime\prime}_{3}y^{\prime})\}+ip^{\prime\prime\prime}_3y\Bigr)=
$$
$$
=-i\{(q_{0}y^{\prime\prime\prime})^{(4)}+
(q_{0}y^{(4)})^{\prime\prime\prime}\}-(p_0y^{\prime\prime\prime})^{\prime\prime\prime}-$$
$$
-i\{(q_{1}y^{(\prime\prime)})^{\prime\prime\prime}+(q_{1}y^{\prime\prime\prime})^{\prime\prime}\}
+(p^{\prime}_1y^{\prime\prime})^{\prime\prime}+i
\{(q^{\prime}_{2}y^{\prime})^{\prime\prime}+(q^{\prime}_{2}y^{\prime\prime})^{\prime  }\}-(p^{\prime\prime}_2y^{\prime})^{\prime}-i
\{(q^{\prime\prime}_{3}y )^{\prime}+(q^{\prime\prime}_{3}y^{\prime})\}+p^{\prime\prime\prime}_3y.
$$

Вернемся к общему случаю.
 Напомним, что, согласно определению (см., например, \cite{EM}),  все квазипроизводные $y^{[j]}$, $j=0,1,...,2n ,$  функции $y$ локально абсолютно непрерывны на $(a,b)$. По условию А) таковой является и функция $q_0$. Поэтому функция $y^{(n)}={y^{[n]}}/{\sqrt{2q_0}}$ (см. формулу (\ref{f10})) является локально абсолютно непрерывной  на $(a,b)$  функцией, и, следовательно,  $y^{(n+1)}\in L^1_{loc}(a,b)$. Учитывая это, первые два слагаемых в правой части формулы (\ref{f201}) преобразуются к виду $(q_0y^{(n)})^{\prime}+q_0y^{(n+1)}-ip_0y^{(n)}(=:\Phi_0)$. Тогда формула (\ref{f201}) примет вид
$$\label{f202}
y^{[n+1]}=
\Phi_0-i\sum\limits_{j=1}^{n}(-1)^j\tilde{\varphi_j}y^{(n-j)}.
$$

Взяв в формуле (\ref{f20}) $k=1$ и учитывая последнюю формулу, после элементарных преобразований получаем, что
$$
y^{[n+2]}=
\Phi_1-i\sum\limits_{j=2}^{n}(-1)^j(\tilde{\varphi_j}y^{(n-j)})^{\prime}-
\sum\limits_{j=1}^{n-1}(-1)^{j+2}C_{j+1}^{1}\left[ip_{j+1}+(1-\frac{2}{j+1})q_{j+1}\right]y^{(n-j)},\,\,
$$
где $$\Phi_1=\Phi^{\prime}_0+(q_1y^{(n-1)})^{\prime}+q_1y^{(n)}+ip^{\prime}_1y^{(n-1)}.$$
Продолжая эти рассуждения, методом математической индукции легко доказать,что
$$
y^{[n+k]}=\Phi_{k-1}-i\sum\limits_{j=k}^{n}(-1)^j(\tilde{\varphi_j}y^{(n-j)})^{(k-1)}-$$
$$-
\sum\limits_{j=k-1}^{n-1} \left \{\sum\limits_{s=1}^{k-1}((-1)^{j+s+1}C_{j+1}^{s}\left[ip_{j+1}+(1-\frac{2s}{j+1})q_{j+1}\right]y^{(n+s-1-j)})^{(k-s-1)} \right\},$$
$$ k=1,2,\ldots,n,
$$
где
 $$\Phi_{0}=(q_0y^{(n)})^{\prime}+q_0y^{(n+1)}-ip_0y^{(n)}$$и
$$
\Phi_{s}=\Phi^{\prime}_{s-1}+(-1)^{s+1}[\{(q^{(s-1)}_{s}y^{(n-s)})^{\prime}+q^{(s-1)}_{s}y^{(n+1-s)}\}+ip^{(s)}_sy^{(n-s)}],\,\,\,s=1, 2,\ldots,n-1,
$$
и поэтому выражение $\tau_{2n+1}$ (см. \eqref{f30}) имеет вид
$$
\tau_{2n+1} y=(-1)^ni
\{(q_{0}y^{(n+1)})^{(n)}+(q_{0}y^{(n)})^{(n+1)}\}+\sum\limits_{k=0}^{n}(-1)^{n+k}(p^{(k)}_ky^{(n-k)})^{(n-k)}+$$
$$+i\sum\limits_{k=1}^{n}(-1)^{n+k+1}\{(q^{(k)}_{k}y^{(n+1-k)})^{(n-k)}+(q^{(k)}_{k}y^{(n-k)})^{(n+1-k)}\}.
$$
Тем самым, дифференциальное выражение с коэффициентами-распределениями вида \eqref{0}  совпадает с регуляризованным дифференциальным выражением \eqref{f30}, с помощью которого можно определять соответствующие  минимальные и максимальные операторы и получать утверждения такого же типа, как в \cite{MSh}.


\medskip
К.\,A.~Mirzoev

Lomonosov Moscow State University,

Department of Mechanics and Mathematics

email: \ mirzoev.karahan@mail.ru

\medskip
A.\,A.\, Shkalikov

Lomonosov Moscow State University,

Department of Mechanics and Mathematics

email: shkalikov@mi.ras.ru


\begin{thebibliography}{99}

\bibitem{MSh}
Мирзоев\,К.\,А.,\,Шкаликов \,А.\,А.\,  Дифференциальные операторы четного порядка с коэффициентами-распределениями // Математические заметки,  2016.\, 99, 5.\,с.788-793. Engl. translation in:
K. A. Mirzoev, A. A. Shkalikov, {\emph Differential Operators of Even Order with Distribution Coefficients}, Math. Notes, 99:5 (2016), 779--784.

\bibitem{EM}
 W.\,N.~Everitt, L.~Marcus
Boundary Value Problems and Sympletic Algebra for Ordinary Differential and Quasi-Differential Operators.
 AMS,\, Mathematical Surveys and Monographs,
 1999, v.61.

\end{thebibliography}
\end{document}